\documentstyle{article}

\newtheorem{theorem}{Theorem}[section]
\newtheorem{lemma}{Lemma}[section]

\begin{document} 

\begin{center} 
\LARGE\bf 
On Storage Operators\\[1cm]
\end{center}

\begin{center} 
\bf 
Karim NOUR\\
\rm
LAMA - Equipe de Logique\\
Universit\'e de Savoie\\
73376 Le Bourget du Lac\\
e-mail nour@univ-savoie.fr\\[1cm]
\end{center}
\begin{abstract} 
In 1990  Krivine (1990b) introduced the notion of storage operators. They are
$\lambda$-terms which simulate call-by-value in the call-by-name strategy.  Krivine (1990b) has
shown that there is a very simple type in the $AF2$ type system for storage operators
using G\H{o}del translation from classical to intuitionistic logic. Parigot (1993a) and  Krivine
(1994) have shown that storage operators play an important tool in classical logic. In this
paper, we present a synthesis of various results on this subject. 
\end{abstract}

\section{Introduction}

Lambda-calculus as such is not a computational model. A reduction strategy is needed. In this paper,
we consider $\lambda$-calculus with the left reduction. This strategy has much advantages : it always
terminates when applied to a normalizable $\lambda$-term and it seems more economic since we compute a
$\lambda$-term only when we need it. But the major drawback of this strategy is that a function must
compute its argument every time it uses it. This is the reason why this strategy is not really
used. In 1990  Krivine (1990b) introduced the notion of storage operators in order to avoid this
problem and to simulate call-by-value when necessary. \\

The $AF2$ type system is a way of interpreting the proof rules for the second order intuitionistic
logic plus equational reasoning as construction rules for terms. Krivine (1990b) has shown that, by
using G\H{o}del translation from classical to intuitionitic logic (denoted by $^g$), we can find in
system $AF2$ a very simple type for storage operators. Historically the type was discovered before the
notion of storage operator itself. Krivine (1990a) proved that as far as totality of functions is
concerned second order classical logic is conservative over second order intuitionistic logic. To
prove this,  Krivine introduced the following notions : $A[x]$ is an input (resp. output) data type
if one can prove intuitionistically $A[x] \rightarrow A^g[x]$ (reps. $A^g[x] \rightarrow \neg \neg
A[x])$. Then if $A[x]$ is an input data type and $B[x]$ is an output data type, then if one prove
$A[x] \rightarrow B[x]$ classically one can prove it  intuitionistically. The notion of storage
operator was discovered by investigating the property of all $\lambda$-terms of type $N^g[x]
\rightarrow \neg \neg N[x]$ where $N[x]$ is the type of integers.\\

Parigot (1992) and  Krivine (1994) have extended the system $AF2$ to the classical logic. The method
of  Krivine is very simple : it consists of adding a new constant, denoted by $C$, with the
declaration $C : \forall X \{ \neg \neg X \rightarrow X \}$ which axiomatizes classical logic over
intuitionistic logic. For the constant $C$, he adds a new reduction rule which is a particular case
of a rule given by Felleisen (1987) for control operator. Parigot considerd a (second order) naturel
deduction system with several conclusions which is more convenient that the usual naturel deduction
system with the classical absurdity rule. Its computational interpretation is a natural extation of
$\lambda$-calculus, called $\lambda \mu$-calculus, which preserves the main properties of
$\lambda$-calculus and alows to model controle structures too. In these systems the property of the
unicity of representation of data is lost, but  Parigot (1993a) and Krivine (1994) have shown that
storage operators typable in $AF2$ can be used to find the values of classical integers.\\

This paper studies some properties of storage operators in pure and typed $\lambda$-calculus. We
present, in particular, the results of Krivine, Parigot and the author.

\section{Pure and typed $\lambda$-calculus}

Let $t,u_1,...,u_n$ be $\lambda$-terms, the application of $t$ to $u_1,...,u_n$ is denoted by 
$(t)u_1...u_n$. $Fv(t)$ is the set of free variables of a $\lambda$-term $t$. The $\beta$-reduction
(resp. $\beta$-equivalence) relation is denoted by $u \rightarrow_{\beta} v$ (resp. $u \simeq_{\beta}
v$). If $t$ is a normalizable $\lambda$-term, we denote by $N(t)$, the number of steps used to go
from $t$ to its normal form. The notation $\sigma(t)$ represents the result of the simultaneous
substitution $\sigma$ to the free variables of $t$ after a suitable renaming of the bounded variables
of $t$. We denote by $(u)^n v$ the $\lambda$-term $(u)...(u)v$ where $u$ occurs $n$ times, and
$\overline{u}$ the sequence of $\lambda$-terms $u_1,...,u_n$ $(n \geq 0)$. If $\overline{u} =
u_1,...,u_n$, we denote by $(t)\overline{u}$ the $\lambda$-term $(t)u_1...u_n$.  \\  Let us recall
that a $\lambda$-term $t$ either has a head redex [i.e. $t=\lambda x_1 ...\lambda x_n (\lambda x u) v
\overline{v}$, the head redex being $(\lambda x u) v$], or is in head normal form [i.e. $t=\lambda
x_1 ...\lambda x_n (x) \overline{v}$]. The notation $u \succ v$ means that $v$ is obtained from $u$
by some head reductions.  A $\lambda$-term t is said to be solvable if and only if the head reduction
of $t$ terminates. If $u \succ v$, we denote by $n(u,v)$ the length of the head reduction between $u$
and $v$. And if $t$ is solvable, we denote by $n(t)$ the number of steps used to go from $t$ to its
head normal form. Krivine (1990b) has shown that :

\begin{lemma}
1) If $u \succ v$, then, for any substitution $\sigma$, $\sigma(u) \succ \sigma(v)$, and $n(\sigma(u),\sigma(v))$=n(u,v). \\ 
2) If $u \succ v$, then, for every sequence of $\lambda$-terms $\overline{w}$, there is a $w$, such that
$(u)\overline{w} \succ w$, $(v)\overline{w} \succ w$, and $n((u)\overline{w},w)=n((v)\overline{w},w)+n(u,v)$.
\end{lemma}

Lemma 2.1 shows that to make the head reduction of $\sigma(u)$ (resp. of $(u)\overline{w}$)
it is equivalent to make some steps in the head
reduction of $u$, and after make the head reduction of $\sigma(v)$ (resp. of $(v)\overline{w}$). \\

The types will be formulas of second order predicate logic over a given
language. The logical connectives are $\perp$ (a predicate symbol 0-air for absurde), $\rightarrow$,
and $\forall$. There are individual (or first order) variables denoted by $x,y,z,...,$ and predicate
(or second order) variables denoted by $X,Y,Z,....$ We do not suppose that the language has a special
constant for equality. Instead, we define the formula $u=v$ (where $u,v$ are terms) to be $\forall
Y(Y(u) \rightarrow Y(v))$ where $Y$ is a unary predicate variable. Such a formula will be called an
equation. We denote by $a \approx b$ the equivalence binary relation such that :  if $a=b$ is an
equation, then $a[t_1 / x_1,..., t_n / x_n] \approx b[t_1 / x_1,...,t_n / x_n]$. The formula $F_1
\rightarrow (F_2 \rightarrow(...\rightarrow (F_n \rightarrow G)...))$ is also denoted by
$F_1,F_2,...,F_n \rightarrow G$. For every formula $A$, we denote by $\neg A$ the formula $A
\rightarrow \perp$.  \\
Let $t$ be a $\lambda$-term, $A$ a type, $\Gamma = x_1 : A_1 ,..., x_n :
A_n$ a context, and $E$ a set of equations. We define by means of the following rules the notion
``$t$ is of type $A$ in $\Gamma$ with respect to $E$'' ; this notion is denoted by 
$\Gamma\vdash_{AF2} t:A$.  

\begin{center} 
(1) $\Gamma\vdash_{AF2} x_i:A_i$ $(1\leq i\leq n)$
\end{center} 

\begin{minipage}[t]{170pt}
$  ( 2 ) \quad \displaystyle\frac{ \Gamma,x:A \vdash_{AF2} t:B } { \Gamma\vdash_{AF2} \lambda xt:A \rightarrow B }$ \\
\end{minipage} 
\begin{minipage}[t]{170pt}\sl 
$  ( 3 ) \quad  \displaystyle\frac{ \Gamma\vdash_{AF2} u:A \rightarrow B \quad \Gamma\vdash_{AF2} v:A} { \Gamma\vdash_{AF2}
(u)v:B }$ \\ \end{minipage}

\begin{minipage}[t]{170pt}
$  ( 4 ) \quad \displaystyle\frac{ \Gamma\vdash_{AF2} t:A } { \Gamma\vdash_{AF2} t:\forall xA }$ {\rm (*)}\\
\end{minipage} 
\begin{minipage}[t]{170pt}\sl 
$  ( 5 ) \quad  \displaystyle\frac{ \Gamma\vdash_{AF2} t:\forall xA} { \Gamma\vdash_{AF2} t:A[u/x] }$  {\rm (**)}\\ 
\end{minipage}

\begin{minipage}[t]{170pt}
$  ( 6 ) \quad \displaystyle\frac{ \Gamma\vdash_{AF2} t:A } { \Gamma\vdash_{AF2} t:\forall XA }$  {\rm (*)}\\
\end{minipage} 
\begin{minipage}[t]{170pt}\sl 
$  ( 7 ) \quad  \displaystyle\frac{ \Gamma\vdash_{AF2} t:\forall XA} { \Gamma\vdash_{AF2} t:A[G/X] }$  {\rm (**)}\\ 
\end{minipage}

\begin{center}
$  ( 8 ) \quad \displaystyle\frac{ \Gamma\vdash_{AF2} t:A[u/x] \quad u \approx v} { \Gamma\vdash_{AF2}
t:A[v/x] }$ \\
\end{center}

With the following conditions : (*) $x$,$X$ have no free occurence in $\Gamma$ and
(**) $u$ (resp. $G$) is a term (resp. formula). \\
This typed $\lambda$-calculus system is called $AF2$ (for {\it Arithm\'etique Fonctionnelle du
second ordre}). It has the following properties (Krivine 1990a).
 
\begin{theorem}
1) Types are preserved during reduction.\\
2) Typable $\lambda$-terms are strongly normalizable.
\end{theorem}

\section{Storage operators}

For every $n \in {\bf N}$, we define the Church integer $\underline{n} = \lambda x
\lambda f (f)^n x$. Let $\underline{s} = \lambda n\lambda x\lambda f((n)(f)x)f$ ; it is easy to check
that $\underline{s}$ is a $\lambda$-term for the successor.\\

Let $F$ be a $\lambda$-term (a function). During the computation, by left reduction, of $(F)\theta_n$ (where
$\theta_n \simeq_{\beta} \underline{n}$), $\theta_n$ may be computed as many times as $F$ uses
it. We would like to transform $(F)\theta_n$ to $(F)\underline{n}$. We also want this transformation
depends only on $\theta_n$ (and not $F$). In other words we look for some closed $\lambda$-terms $T$
with the following properties : \\
- For every $\lambda$-term $F$, $n \in {\bf N}$, and $\theta_n \simeq_{\beta} \underline{n}$, we have
$(T)\theta_nF \succ (F)\underline{n}$; \\ 
- The computation time of $(T)\theta_nF \succ (F)\underline{n}$ depends only on $\theta_n$.  \\

{\bf Definition (temporary)} : A closed $\lambda$-term $T$ is called storage operator for Church
integers iff for every $n \in {\bf N}$, and for every $\theta_n \simeq_{\beta}
\underline{n}$, $(T)\theta_n f \succ (f)\underline{n}$ (where $f$ is a new variable). \\

It is clear that a storage operator satisfies the required properties. Indeed, since we have
$(T)\theta_n f \succ (f)\underline{n}$, then the variable $f$ never comes in head position during the
reduction, and we may then replace $f$ by any $\lambda$-term. We will show (see Theorem 3.1) that it
is not possible to get the normal form of $\theta_n$. We then change the definition.\\

{\bf Definition (temporary)} :	A closed $\lambda$-term $T$ is called storage operator for Church
integers iff for every  $n \in {\bf N}$, there is a closed $\lambda$-term
$\tau_n\simeq_{\beta} \underline{n}$, such that for every $\theta_n \simeq_{\beta}
\underline{n}$, $(T)\theta_n f \succ (f)\tau_n$ (where $f$ is a new variable). \\ 

Krivine (1990b) has shown that, by using G\H{o}del translation from classical to intuitionitic
logic, we can find a very simple type for storage operators. But the $\lambda$-term $\tau_n$
obtained may contain variables substituted by $\lambda$-terms $u_1,...,u_m$ depending on
$\theta_n$. Since the $\lambda$-term $\tau_n$ is $\beta$-equivalent to $\underline{n}$,
therefore, the left reduction of the $\tau_n[u_1/x_1,...,u_m/x_m]$ is equivalent to the left
reduction of $\tau_n$ and the $\lambda$-terms $u_1,...,u_m$ will therefore never be evaluated
during the reduction. \\

{\bf Definition (final)} : A closed $\lambda$-term $T$ is called a storage operator for Church
integers iff for every $n \in {\bf N}$, there is a $\lambda$-term $\tau_n\simeq_{\beta}
\underline{n}$, such that for every $\theta_n \simeq_{\beta} \underline{n}$, there is a substitution
$\sigma$, such that $(T)\theta_nf \succ (f)\sigma(\tau_n)$ (where $f$ is a new variable). \\

Let $F$ be any $\lambda$-term (for a function), and $\theta_n$ a $\lambda$-term $\beta$-equivalent to
$\underline{n}$. During the computation of $(F)\theta_n$, $\theta_n$ may be computed each time it comes in head
position. Instead of computing $(F)\theta_n$, let us look at the head reduction of $(T)\theta_n F$. Since
it is $\{(T)\theta_n f\}[F/f]$, by Lemma 2.1, we shall first reduce $(T)\theta_n f$ to its head normal
form, which is $(f)\sigma(\tau_n)$, and then compute $(F)\sigma'(\tau_n)$ ($\sigma'=\gamma \circ \sigma$
where $\gamma (f)=F$ an d$\gamma (x)=x$ if $x \not = f$). The computation has been decomposed into
two parts, the first being independent of $F$. This first part is essentially a computation of
$\theta_n$, the result being $\tau_n$, which is a kind of normal form of $\theta_n$. The substitutions made
in $\tau_n$ have no computational significance, since $\underline{n}$ is closed. So, in the computation of
$(T)\theta_n F$, $\theta_n$ is computed first, and the result is given to $F$ as an argument, $T$ has
stored the result, before giving it, as many times as needed, to any function. \\

If we take : $T_1 = \lambda n((n)\delta)G$ where
$\delta = \lambda f(f)\underline{0}$ and $G = \lambda x\lambda y(x)\lambda z(y)(\underline{s})z$ ; \\
$T_2 = \lambda n\lambda f(((n)f)F)\underline{0}$ where $F = \lambda x\lambda y(x)(\underline{s})y$, 
then we can check that for every $\theta_n \simeq_{\beta} \underline{n}$, $(T_i)\theta_n f \succ
(f)(\underline{s})^n \underline{0}$ ($i=1$ or $2$) (Krivine 1990a and Nour 1993a). Therefore $T_1$
and $T_2$ are storage operators for Church integers.\\

The most effective  storage operators for Church integers - found by 
Krivine  - give as result $(\underline{s})^n\underline{0}$. A question arises : {\it Can we find
storage operators for Church integers which give normal forms as result ?} This kind of storage
operators are called strong storage operators. We have shown (Nour 1995a) that :
 
\begin{theorem} Church integers do not have strong storage operators.
\end{theorem}

The nonexistence of strong storage operators for Church integers results from the following facts:\\
- {\it The infinity of integers} : We can prove that every finite subset of Church integers has
strong storage operators (Nour 1995a). \\
- {\it The representation of integers} : We can prove that we cannot create a Church integer
$\underline{n}$ ($n \geq 1$) during head reduction in the application. If we change the representation
of integers, we can find strong storage operators. For every $n \in {\bf N}$, we define the recursive
integer $\overline{n}$ by induction : $\overline{0}=\lambda f\lambda xx$ and $\overline{n+1}=\lambda
f\lambda x(f)\overline{n}$. Let $\overline{s}=\lambda n\lambda f\lambda x(f)n$ ; it is easy to check
that $\overline{s}$ is a $\lambda$-term for successor. If we take $T'=\lambda \nu(\nu)\rho \tau \rho$
where $\tau=\lambda f(f) \overline{0}$, $\rho=\lambda y\lambda z(G)(y)z \tau z$, and $G=\lambda
x\lambda y(x)\lambda z(y) \lambda f\lambda x(f)z$, then, for every $\theta_n \simeq_{\beta}
\overline{n}$, $(T')\theta_nf \succ (f)\overline{n}$. Therefore $T'$ is a strong storage operators
for recursive integers (Nour 1995a).

\section{Directed $\lambda$-calculus and storage operators}

A closed $\lambda$-term $T$ is a storage operator for Church integers iff for every $n \in
{\bf N}$, there is a $\lambda$-term $\tau_n \simeq_{\beta} \underline{n}$, such that for every $\theta_n
\simeq_{\beta} \underline{n}$, there is a substitution $\sigma$, such thatÊ$(T)\theta_n f \succ
(f)\sigma(\tau_n)$. Let's analyse the head reduction $(T)\theta_n f \succ (f)\sigma(\tau_n)$, by
replacing each $\lambda$-term which comes from $\theta_n$ by a new variable. This will help us to
better understand the  Krivine proof of his principal storage Theorem (Theorem 5.2)
and also to justify the introduction of directed $\lambda$-calculus which allows to find  similar
results in the general case.\\
  
If  $\theta_n \simeq_{\beta} \underline{n}$,  then $\theta_n \succ \lambda x \lambda g(g)t_{n-1}$,  $t_{n-k} \succ (g)t_{n-k-1}$ 
$(1 \leq k \leq n-1)$,  $t_0 \succ x$,  and  $t_k \simeq_{\beta} (g)^k x$ $(0 \leq k \leq n-1)$. Let
$x_n$ be a new variable ($x_n$ represents $\theta_n$). $(T)x_nf$ is solvable, and its head normal form
does not begin by $\lambda$, therefore it is a variable applied to some arguments. The free variables of
$(T)x_n f$ are $x_n$ and $f$, we then have two possibilities for its head normal form : $(f)\delta$
(in this case we stop) or $(x_n)a_1...a_m$. Assume we obtain $(x_n)a_1...a_m$. The variable $x_n$
represents $\theta_n$, and  $\theta_n \succ \lambda x \lambda g(g)t_{n-1}$, therefore $(\theta_n)a_1...a_m$ and
$((a_2)t_{n-1}[a_1/x,a_2/g])a_3...a_m$ have the same head normal form. The $\lambda$-term
$t_{n-1}[a_1/x,a_2/g]$ comes from $\theta_n$. Let $x_{n-1,a_1,a_2}$ be a new variable
$(x_{n-1,a_1,a_2}$ represents $t_{n-1}[a_1/x,a_2/g])$. The $\lambda$-term
$((a_2)x_{n-1,a_1,a_2})a_3...a_m$ is solvable, and its head normal form does not begin by $\lambda$,
therefore it is a variable applied to some arguments. The free variables of
$((a_2)x_{n-1,a_1,a_2})a_3...a_m$ are among $x_{n-1,a_1,a_2}$, $x_n$, and $f$, we then have three
possibilities for its head normal form :  $(f)\delta$ (in this case we stop) or $(x_n)b_1...b_r$ or
$(x_{n-1,a_1,a_2})b_1...b_r$. Assume we obtain $(x_{n-1,a_1,a_2})b_1...b_r$. The
variable$(x_{n-1,a_1,a_2}$ represents $t_{n-1}[a_1/x,a_2/g]$, and $t_{n-1} \succ (g)t_{n-2}$,
therefore $(t_{n-1}[a_1/x,a_2/g])b_1...b_r$ and $((a_2)t_{n-2}[a_1/x,a_2/g])b_1...b_r$ have the
same head normal form. The $\lambda$-term $t_{n-2}[a_1/x,a_2/g]$ comes from $\theta_n$. Let
$x_{n-2,a_1,a_2}$ be a new variable $(x_{n-2,a_1,a_2}$ represents $t_{n-2}[a_1/x,a_2/g])$. The
$\lambda$-term $((a_2)x_{n-2,a_1,a_2})b_1...b_r$ is solvable, and its head normal form does not begin by
$\lambda$, therefore it is a variable applied to arguments. The free variables of
$((a_2)x_{n-2,a_1,a_2})b_1...b_r$ are among $x_{n-2,a_1,a_2}$, $x_{n-1,a_1,a_2}$, $x_n$, and $f$,
therefore we have four possibilities for its head normal form : $(f)\delta$ (in this case we stop)
or $(x_n)c_1...c_s$ or $(x_{n-1,a_1,a_2})c_1...c_s$ or $(x_{n-2,a_1,a_2})c_1...c_s$ ... and so on...
Assume we obtain $(x_{0,d_1,d_2})e_1...e_k$ during the construction. The variable $x_{0,d_1,d_2}$
represents $t_0[d_1/x,d_2/g]$, and $t_0 \succ x$, therefore $(t_0[d_1/x,d_2/g])e_1...e_k$ and
$(d_1)e_1...e_k$ have the same head normal form ; we then follow the construction with the
$\lambda$-term $(d_1)e_1...e_k$. The $\lambda$-term $(T)\theta_n f$ is solvable, and has $(f)\sigma(\tau)$ as head
normal form, so this construction always stops on $(f)\delta$. We can prove by a simple argument
that $\delta \simeq_{\beta} \underline{n}$.  \\

According to the previous construction, the reduction $(T)\theta_nf \succ (f)\sigma(\tau_n)$ can
be divided into two parts : a reduction that does not depend on $n$ and a reduction that depends on
$n$ (and not on $\theta_n$). If we allow some new reduction rules to get the later reductions,
(something as : $(x_n)a_1a_2 \succ (a_2)x_{n-1,a_1,a_2}$ ; $x_{i+1,a_1,a_2} \succ (a_2)u_{i,a_1,a_2}$
($i>0$) ; $x_{0,a_1,a_2} \succ  a_1$) we obtain an equivalent definition for the storage
operators for Church integers : a closed $\lambda$-term $T$ is a storage operator for Church integers
iff for every $n \in {\bf N}$, $(T)x_nf \succ (f)\delta_n$ where $\delta_n \simeq_{\beta}
\underline{n}$. To prove his storage Theorem (Theorem 5.2),  Krivine used the sufficient condition of
the laste equivalence.\\

The notion of storage operators can be generalized for each set of closed normal $\lambda$-terms.\\

Let $t$ be a closed normal $\lambda$-term and $T$ a closed $\lambda$-term. We sad that
$T$ is a storage operator for $t$ iff there is a $\lambda$-term $\tau_t \simeq_{\beta}
t$, such that for every  $\lambda$-term $\theta_t \simeq_{\beta} t$, there is a substitution $\sigma$, such
that $(T)\theta_t f \succ (f)\sigma(\tau_t)$ (where $f$ is a new variable).
Let $D$ be set of closed normal $\lambda$-terms and $T$ a closed $\lambda$-term. We sad that $T$ is
a storage operator for $D$ iff it is a storage operator for every $t$ in $D$.\\

The directed $\lambda$-calculus is an extension of the ordinary $\lambda$-calculus built for tracing a normal
$\lambda$-term $t$ during some head reduction. Assume $u$ is some normal $\lambda$-term having $t$ as a
subterm. We wish to trace the places where we really have to know what $t$ is during the reduction
of $u$. We will present how the directed $\lambda$-calculus allows to find an equivalent -and easily
expressed - definition for the storage operators. \\

Let $V$ be a set of variables of pure $\lambda$-calculus. The set of terms of
directed $\lambda$-calculus, denoted by $\Lambda[]$, is defined in the following way : \\
- If $x \in V$, then $x \in \Lambda[]$ ;  \\ 
- If $x \in V$, and $u\in \Lambda[]$, then $\lambda x u \in \Lambda[]$ ;  \\
 - If $u,v \in \Lambda[]$, then $(u)v \in \Lambda[]$ ;  \\
- If $t \in \Lambda$ is a normal $\lambda$-term, such that $Fv(t) \subseteq \{x_1,...,x_n \}$, and
$a_1,...,a_n \in \Lambda[]$, then $[t]<a_1/x_1,...,a_n/x_n> \in \Lambda[]$.\\
A $\lambda[]$-term of the form $[t]<a_1/x_1,...,a_n/x_n>$ is said to be a  box directed by $t$. 
This notation represents, intuitively, the $\lambda$-term $t$ where all free variables $x_1,...,x_n$
will be replaced by $a_1,..,a_n$. The substitution $<a_1/x_1,...,a_n/x_n>$ is denoted by $<{\bf
a}/{\bf x}>$.  \\
A $\lambda[]$-term of the form $(\lambda xu)v$ is called $\beta$-redex ; $u[v/x]$ is
called its contractum. A $\lambda[]$-term of the form $[t]<{\bf a}/{\bf x}>$ is called $[]$-redex
; its contractum $R$ is defined by induction on $t$ :\\
- If $t=x_i$ $(1 \leq i \leq n)$, then $R=a_i$ ;  \\ 
- If $t=x \not = x_i$ $(1 \leq i \leq n)$, then $R=x$ ; \\
- If $t=\lambda xu$, then $R=\lambda y[u]<{\bf a}/{\bf x},y/x>$ where $y \not \in Fv({\bf a})$ ; \\
- If $t=(u)v$, then $R=([u]<{\bf a}/{\bf x}>)[v]<{\bf a}/{\bf x}>$. \\

By interpreting the box $[t]<a_1/x_1,...,a_n/x_n>$ by $t[[ a_1/x_1,...,a_n/x_n ]]$ (the
$\lambda$-term $t$ with an explicit substitution), the new reduction rules are those that allow to
really do the substitution. This kind of $\lambda$-calculus has been studied by Curien (1988) ; his
$\lambda \sigma$-calculus contain terms and substitutions and is intended to better control the
substitution process created by $\beta$-reduction, and then the implementation of the
$\lambda$-calculus. The main difference between the $\lambda \sigma$-calculus and the directed
$\lambda$-calculus is : The first one produces an explicit substitution after each
$\beta$-reduction.  The second only ``executes'' the substitutions given in advance.  We can
therefore consider the directed $\lambda$-calculus as a restriction (the interdiction of producing
explicit substitutions) of $\lambda \sigma$-calculus ; a well adapted way to the study of the head
reduction.\\

Every $\lambda []$-term $t$ can be - uniquely - written as $\lambda x_1...\lambda
x_n(R)t_1...t_m$ $n,m \geq 0$, $R$ being a variable or a redex. If $R$ is a variable, we say that
$t$ is a $\beta []$-head normal form. If $R$ is a redex, we say that $R$ is the head redex
of $t$. The notation $u \succ_{\beta []} v$ means that $v$ is obtained from $u$ by some head
reductions.\\

Now, we can state the Theorem which gives an equivalent definition for storage operators (Nour and
David 1995). 

\begin{theorem}
Let $t$ be a closed normal $\lambda$-term, and $T$ a closed
$\lambda$-term. $T$ is a storage operator for $t$ iff there is a $\lambda$-term $\tau_t
\simeq_{\beta} t$, such that \\
$(T)[t]f \succ_{\beta []} (f)\tau_t[[t_1]<{\bf a_1}/{\bf x_1}>/y_1,...,[t_m]<{\bf a_m}/{\bf x_m}>/y_m]$. 
\end{theorem} 

To prove the necessary condition we associate to every $\theta_t \simeq_{\beta} t$ a special
substitution $S_{\theta}$ over the boxes directed by subterms of $t$ such that
$S_{\theta}([t])=\theta_t$  and satisfying the following property : if $u \succ_{\beta []} v$ then
$S_{\theta}(u) \succ S_{\theta}(v)$. Then $(T)\theta_t f \succ (f)\sigma(\tau_t)$. For the sufficient
condition we use the idea given at the begining of this paragraph. The only difficulty is to prove
that  $\tau_t \simeq_{\beta} t$. For that we use the fact that $\tau_t$ dos not depend on
$\theta_t$.\\

The laste result allows to find some important properties for storage operators (Nour and
David 1995).  

\begin{theorem} 1) Let $D$ be a set of closed normal $\lambda$-terms, $T$ and $T'$ two closed
$\lambda$-terms. If $T$ is a storage operator for $D$, and $T' \simeq_{\beta} T$, then $T'$ also is
a storage operator for $D$. \\
2) The set of  storage operators for a set of closed normal $\lambda$-terms is not
recursive. But the set of  storage operators for a finite set of closed normal $\lambda$-terms is
recursively enumerable.  \\
3) Each finite set of normal $\lambda$-terms having all distinct $\beta\eta$-normal forms has a
storage operator.  \\ 
4) Let $t$ be a closed normal $\lambda$-term, and $T$ a closed $\lambda$-term. If T is
a  storage operator for $t$, then there are two constants $A_{T,t}$ and $B_{T,t}$, such that for
every $\theta_t \simeq_{\beta} t$, $n((T)\theta_tf) \leq A_{T,t} N(\theta_t) + B_{T,t}$. 
\end{theorem}

\section{Storage operators in typed $\lambda$-calculus}

Each data type generated by free algebras can be defined by a second order
formula. The type of integers is the formula : $N[x]= \forall X \{ X(0), \forall y(X(y) \rightarrow
X(sy)) \rightarrow X(x) \}$ where $X$ is a unary predicate variable, $0$ is a constant symbol for zero, and
$s$ is a unary function symbol for successor. The formula $N[x]$ means semantically that $x$ is an
integer iff $x$ belongs to each set $X$ containing $0$ and closed under the successor
function $s$. It is easy to check that, for every $n \in {\bf N}$, the Church integer $\underline{n}$
is of type $N[s^n(0)]$ and $\underline{s}$ is of type $\forall y(N[y] \rightarrow N[sy])$.    \\
A set of equations $E$ is said to be adequate with the type of integers iff :
$s(a) \not \approx 0$ and if $s(a) \approx s(b)$, then $a \approx b$. In the rest of the paper, we
assume that all sets of equations are adequate with the type of integers.\\

The system $AF2$ has the property of the unicity of integers representation (Krivine 1990a). 

\begin{theorem} 
 Let $n \in {\bf N}$, if $\vdash_{AF2} t :N[s^n (0)]$, then $t \simeq_{\beta} \underline{n}$. 
\end{theorem} 

A very important property of data type is the following (we express it for the type
of integers) : in order to get a program for a function $f : {\bf N} \rightarrow {\bf N}$ it is sufficient to prove $\vdash
\forall x ( N[x] \rightarrow N[f(x)] )$. For example a proof of $\vdash \forall x ( N[x] \rightarrow N[p(x)] )$ from the equations
$p(0)=0$, $p(s(x))=x$ gives a $\lambda$-term for the predecessor in Church intergers (Krivine 1990a). \\

If we try to type a storage operator $T$ for Church integers in $AF2$ type system, we naturally
find the type $\forall x \{ N[x] \rightarrow \neg \neg N[x] \}$. But this type does not characterize
the storage operators (take for example $T=\lambda \nu \lambda f(f)\nu$). This comes from the fact
that the type $\forall x \{ N[x] \rightarrow \neg \neg N[x] \}$ does not take into account the
independency of $\tau_n$ from $\theta_n$. To solve this problem, we must prevent the use of the first
$N[x]$ in $\forall x \{ N[x] \rightarrow \neg \neg N[x] \}$ as well as his subtypes to prove the
second $N[x]$. \\

For each formula $F$ of $AF2$, we indicate by $F^g$ the formula obtained by
putting $\neg$ in front of each atomic formulas of $F$ ($F^g$ is called the G\H{o}del translation
of $F$). For example : $N^g[x]=\forall X \{\neg X(0),\forall y(\neg X(y) \rightarrow \neg X(sy)) \rightarrow \neg X(x)
\}$. It is well known that, if $F$ is provable in classical logic, then $F^g$ is provable in
intuitionistic logic (Krivine 1990a).\\ 

We can check that $\vdash_{AF2} T_1,T_2 : \forall x \{N^g[x] \rightarrow \neg\neg N[x] \}$. And, in
general, we have the following Theorem (Krivine 1990a, Nour 1994) :

\begin{theorem} 
If $\vdash_{AF2} T: \forall x\{N^g[x] \rightarrow \neg\neg N[x]\}$, then $T$ is a storage operator for
Church integers. 
\end{theorem}

We will give some ideas for the proofs of this Theorem. Krivine (1990a) introduced a semantic for his
system and he proved that : if $t$ is of type $A$ then $t$ belongs $A$. Since $T$ is of type $\forall
x \{ N^g[x] \rightarrow \neg \neg N[x] \}$  then $T$ belongs $N^g[s^n(0)] \rightarrow \neg \neg
N[s^n(0)]$. With the proper semantic interpretation of $\perp$ we check that $x_n$ belongs 
$N^g[s^n(0)]$ and $f$ belongs $\neg N[s^n(0)]$. This implies that $(T)x_n f$ belongs to $\perp$ which
gives the theorem directly from the choice of the interpretation of $\perp$. We presented (Nour
1994) a syntactical proof of this result. We prove by using only the syntactical properties of the
system $AF2$ that the $\lambda$-term $T$ satisfies the properties which we need.\\

The storage operators given in this paper up to now give as results closed $\lambda$-terms. This kind
of storage operators is called proper storage operators. A question arises : {\it Can we find a typed
non proper storage operator for Church integers ?} We have shown that (Nour 1993b) :

\begin{theorem}
There is a non proper storage operator for Church integers $T$ such that  $\vdash_{AF2} T : \forall x
\{N^g[x] \rightarrow \neg\neg N[x] \}$.
\end{theorem}

An example of a such operator is the following : $T=\lambda \nu(\nu)\gamma D$ where \\
$D=\lambda u\lambda v(u)\lambda w(((\nu)\lambda y(((y)w)u)v)\lambda xx )\lambda g\lambda k\lambda
l(l)\lambda n\lambda m(n)((g)n)m$, \\ $\gamma=\lambda f(((\nu)\lambda
x(f)((((x)n)f)\underline{0})\lambda xx)\lambda x\lambda y\lambda zz$.
                                     
\section{Generalization}

Some authors have been interested in the research of a most general type for storage operators. For
example, Danos and Regnier (1992) have given as type for storage operators the formula $\forall x
\{N^e [x] \rightarrow \neg\neg N[x] \}$ where the operation $e$ is an elaborate G\H{o}del translation
which associates to every formula $F$ the formula $F^e$ obtained by replacing in $F$ each atomic
formula $X(\overline{t})$ by $X_1(\overline{t}),...,X_r(\overline{t}) \rightarrow \perp$. 
Krivine (1993) and the author (Nour 1996a) have given a more general type for storage operators the
formula $\forall x \{ N^G [x] \rightarrow \neg\neg N[x] \}$ where the operation $G$ is the general
G\H{o}del translation which associates to every formula $F$ the formula $F^G$ obtained by replacing
in $F$ each atomic formula $X(\overline{t})$ by a formula $G_X [ \overline{t} / \overline{x} ]$
ending with $\perp$. With the types cited before, we cannot type the simple storage
operator :  $T=\lambda \nu \lambda f(( \nu )\lambda xx)(T_i) \nu f $  ($i=1$ or $2$). This is due
to the fact that the normal form of $T$ contains a variable $\nu$ applied to two arguments
and another $\nu$ applied to three arguments. Therefore, we cannot type $T$  because the
variable $\nu$ is assigned by $N^g[x]$ (for example) and thus the number of the $\nu$-arguments is
fixed once for all. To solve the problem, we replace $N^g[x]$ in the type of storage operators by
another type $N^{\perp}[x]$ which does not limit the number of $\nu$-arguments and only enables to
generate formulas ending with $\perp$ in order to find a general specification for storage
operators.\\

We assume that for every integer $n$, there is a countable set of special
$n$-ary second order variables denoted by $X_{\perp},Y_{\perp},Z_{\perp}$...., and called
$\perp$-variables. A type $A$ is called an $\perp$-type iff $A$ is obtained by the following
rules :   \\
- $\perp$ is an $\perp$-type ;\\
- $X_{\perp}(t_1,...,t_n)$ is an $\perp$-type ;\\
- If $B$ is an $\perp$-type, then $A \rightarrow B$ is an $\perp$-type for every type $A$ ;\\
- If $A$ is an $\perp$-type, then $\forall vA$ is an $\perp$-type for every variable $v$.\\
We add to the $AF2$ type system the new following rules :\\

\begin{minipage}[t]{170pt}
$(6') \quad \displaystyle\frac{\Gamma\vdash t:A } { \Gamma\vdash
t:\forall X_{\perp}A }$ (*)\\
\end{minipage} 
\begin{minipage}[t]{170pt}\sl 
$(7') \quad \displaystyle\frac{\Gamma\vdash t:\forall X_{\perp}A } { \Gamma\vdash t:A[G/X_{\perp}]}$ (**)\\
\end{minipage}

With the following conditions : (*) $X_{\perp}$ has no free occurence in $\Gamma$ and  
(**) $G$ is an $\perp$-type. \\
We call $AF2_{\perp}$ the new type system, and we write $\Gamma\vdash_{\perp} t:A$ if $t$ is typable
in $AF2_{\perp}$ of type $A$ in the context $\Gamma$.\\
We define two sets of types of $AF2$ type system: $\Omega^+$ (set of $\forall$-positive types),
and $\Omega^-$ (set of $\forall$-negative types) in the following way :\\
- If $A$ is an atomic type, then $A \in \Omega^+$, and $A \in \Omega^-$ ;\\
- If $T \in \Omega^+$, and $T' \in \Omega^-$, then, $T' \rightarrow T \in \Omega^+$, and $T \rightarrow T' \in \Omega^-$ ; \\
- If $T \in \Omega^+$ (resp. $T \in \Omega^-$), then $\forall x T \in \Omega^+$ (resp. $\forall x T \in \Omega^-$);\\
- If $T \in \Omega^+$, then $\forall X T \in \Omega^+$ ;\\
- If $T \in \Omega^-$, and $X$ has no free occurence in $T$, then $\forall X T \in \Omega^-$.\\
Therefore, $T$ is a $\forall$-positive types iff the universal second order quantifier appears positively
in $T$. \\
For each predicate variable $X$, we associate an $\perp$- variable $X_{\perp}$. For each
formula $A$ of $AF2$ type system, we define the formula $A^{\perp}$ as follows :  \\
- If $A=R(t_1,...,t_n)$, where $R$ is an $n$-ary predicate symbol, then  $A^{\perp} = A$ ; \\
- If $A=X(t_1,...,t_n)$, where $X$ is an $n$-ary predicate variable, then $A^{\perp} =
X_{\perp}(t_1,...,t_n)$;  \\ 
- If $A=B \rightarrow C$, then $A^{\perp} = B^{\perp} \rightarrow C^{\perp}$ ;\\
- If $A=\forall x B$, then $A^{\perp} = \forall x  B^{\perp}$ ;\\
- If $A=\forall X B$, then $A^{\perp} = \forall X_{\perp} B^{\perp}$.\\
Let $T$ be a closed $\lambda$-term, and $D,E$ two closed types of $AF2$ type system. We say that
$T$ is a storage operator for the pair of types $(D,E)$ iff for every $\lambda$-term
$\vdash_{AF2} t:D$, there are $\lambda$-terms $\tau_t$ and $\tau'_t$, such that $\tau'_t \simeq_{\beta} \tau_t$,
$\vdash_{AF2}\tau'_t:E$, and for every $\theta_t \simeq_{\beta} t$, there is a substitution $\sigma$, such that
$(T)\theta_t f \succ (f)\sigma(\tau_t)$ (where $f$ is a new variable). \\

We have the following generalization (Nour 1995d).

\begin{theorem}
Let $D,E$ be two $\forall$-positive closed types of $AF2$ type system,
such that $E$ does not contain $\perp$. If $\vdash_{\perp} T: D^{\perp} \rightarrow \neg\neg E$, then $T$ is a
storage operator for the pair $(D,E)$.   
\end{theorem}

The condition ``$D,E$ are $\forall$-positive types'' is necessary in order to obtain Theorem 6.1. Indeed,
let $D=\forall X\{ \forall Y(Y \rightarrow X) \rightarrow X \}$, $t=\lambda x(x)\lambda yy$, and $T=\lambda \nu(\nu)\lambda x\lambda f(f)\lambda
y(y)x$. It is easy to check that $D$ is not a $\forall$-positive type, $\vdash_{AF2} t:D$,
$\vdash_{\perp} T:D^{\perp} \rightarrow \neg \neg D$, and $T$ is not a
storage operator for $D$ (Nour 1993a). This counter example also works with the original G\H{o}del
translation and with any general G\H{o}del translation.\\

Theorem 6.1 allows also to generalize the result of  Krivine (Theorem 5.2) to every data type
(booleans, lists, trees, product and sum of data types, ...).

\section{Pure and typed $\lambda C$-calculus}

We add a constant $C$ to the pure
$\lambda$-calculus and we denote by $\lambda C$ the set of new terms also called $\lambda C$-terms.
We consider the following rules of reduction, called rules of head $C$-reduction.\\
(1) $(\lambda  x u) t t_1 ... t_n \rightarrow (u[t / x]) t_1 ... t_n$ for every $u, t, t_1,...,t_n \in
\Lambda  C$. \\
(2) $(C) t t_1 ... t_n \rightarrow (t) \lambda  x (x)t_1 ... t_n$ for every $ t, t_1,...,t_n
\in \Lambda  C$, $x$ being a $\lambda$-variable not appearing in $t_1,...,t_n$.\\
The rule $(2)$ is a particular case of a general law of reduction for control operators given in
(Felleisein 1987) which is $E[Ct/x] \rightarrow (t)\lambda xE$.\\
For any $\lambda C$-terms $t,t'$, we shall write $t \succ_C t'$ if $t'$ is obtained from $t$
by applying these rules finitely many times. \\
A $\lambda C$-term $t$ is said to be $\beta$-normal iff $t$ does not contain a
$\beta$-redex.\\
A $\lambda C$-term $t$ is said to be  $C$-solvable iff  $t \succ_C (f)t_1,...,t_n$
where $f$ is a variable.\\
We add to the $AF2$ type system the new following rule : 
\begin{center} 
(0) $\Gamma \vdash C : \forall X \{ \neg \neg X \rightarrow X \}$
\end{center}
This rule axiomatizes the classical over the intuitionistic logic.
We call $C2$ the new type system, and we write $\Gamma \vdash_{C2} t : A$ if $t$ is of type $A$ in the
context $\Gamma$. In this system we have only the following weak properties (Krivine 1994).

\begin{theorem} 
1) If $\Gamma \vdash_{C2} t:A$, and $t \rightarrow_{\beta} t'$, then $\Gamma \vdash_{C2} t':A$.\\
2) If $\Gamma \vdash_{C2} t:\perp$, and $t \succ_C t'$, then $\Gamma \vdash_{C2} t':\perp$.\\
3) If $A$ is an atomic type, and $\Gamma \vdash_{C2} t:A$, then $t$ is $C$-solvable.
\end{theorem}

In this system, the problem is : given a typed term in classical logic, what kind of program is it
? We shall take the example of integers. Let us call a $\lambda C$-term $\theta$ a classical integer if
$\vdash_{C2}\theta : N[s^n0]$. If $\vdash_{AF2}\theta : N[s^n0]$, then we know that $\theta \simeq_{\beta} \underline{n}$,
and thus we know the operational behaviour of $\theta$. But when $\theta$ is a classical integers, it
is no longer true that $\theta \simeq_{\beta} \underline{n}$. For example $\vdash_{C2} \theta_1=\lambda x \lambda  f(C)\lambda 
y(y)(f)(C)\lambda  z(y)(f)x : N[s0]$. In order to recognize the integer $n$ hidden inside $\theta$
(the value of $\theta$), we have make use of storage operators. Krivine (1994) has shown that :

\begin{theorem}
If $\vdash_{AF2} T: \forall x\{N^g[x] \rightarrow \neg\neg N[x]\}$, then for every $n \in {\bf N}$, there
is a $\lambda $-term $\tau_n \simeq_{\beta} \underline{n}$ such that for every classical integer
$\theta_n$ of value $n$, there is a substitution $\sigma$ such that $(T)\theta_n f \succ_C
(f)\sigma(\tau_n)$ (then $(T)\theta_n \lambda  xx \succ_C \sigma'(\tau_n) \rightarrow_{\beta}
\underline{n}$).  
\end{theorem}

The difficulties to prove this theorem (by comparasion to the Theorem 5.1) are : the operational
characterization of classical integers and the fact that this characterization corresponds to the
behavior of typed storage operators.\\

Theorem 7.2 cannot be generalized for the system $C2$. Indeed, let
$T=\lambda  \nu \lambda  f (f) (C)(T_i)\nu$ ($i=1$ or $2$).  We have
$\vdash_{C2} T: \forall x\{N^g[x] \rightarrow \neg\neg N[x]\}$ and
there is not a $\lambda C$-term $\tau_n \simeq_{\beta} \underline{n}$ such that for every classical
integer $\theta_n$ of value $n$, there is a substitution $\sigma$, such that $(T)\theta_n f \succ_C
(f)\sigma(\tau_n)$ (Nour 1997a). \\

The Theorem 7.2 suggests many questions :\\
- What is the relation between classical integers and the
type $N^g[x]$ ? \\
- Why do we need intuitionistic logic to modelize the storage operators and
classical logic to modelize the control operators ? 

\section{The $M2$ type system}

In this section, we present a new classical type system based on a logical system called
mixed logic. This system allows essentially to distinguish between classical proofs and
intuitionistic proofs. We assume that for every integer $n$, there is a countable set of special
$n$-ary second order variables denoted by $X_C,Y_C,Z_C$...., and called classical variables.\\ 

Let $X$ be an $n$-ary predicate variable or predicate symbol. A type $A$ is said to be
 ending with $X$ iff $A$ is obtained by the following rules :\\
- $X(t_1,...,t_n)$ ends with $X$;  \\ 
- If $B$ ends with $X$, then $A \rightarrow B$ ends with $X$ for every type $A$ ; \\ 
- If $A$ ends with $X$, then $\forall vA$ ends with $X$ for every variable $v$. \\
A type $A$ is said to be a classical type iff $A$ ends with $\perp$ or a
classical variable.\\
We add to the $AF2$ type system the new following rules :

\begin{center}
$(0') \quad \Gamma\vdash C : \forall X_C \{ \neg \neg X_C \rightarrow X_C \}$ \\
\end{center}

\begin{minipage}[t]{170pt}
$  ( 6'' ) \quad \displaystyle\frac{\Gamma\vdash t:A } { \Gamma\vdash t:\forall
X_C A }$ (*)\\
\end{minipage} 
\begin{minipage}[t]{170pt}\sl
$  ( 7'' ) \quad \displaystyle\frac{\Gamma\vdash t:\forall X_C A } { \Gamma\vdash t:A[G/X_C]}$ (**)\\
\end{minipage}

With the following conditions : (*) $X_C$ has no free occurence in $\Gamma$ and 
(**) $G$ is a classical type. \\
We call $M2$ the new type system, and we write $\Gamma \vdash_{M2} t:A$ if $t$ is of type $A$ in the
context $\Gamma$.

\subsection{Properties of $M2$}

With each classical variable $X_C$, we associate a special variable $X^{\bullet}$ of
$AF2$ having the same arity as $X_C$. For each formula $A$ of $M2$, we define the formula
$A$* of $AF2$ in the following way :\\
 - If $A=D(t_1,...,t_n)$ where $D$ is a predicate symbol or a predicate variable, then $A$*=$A$
;\\
- If $A=X_C(t_1,...,t_n)$, then $A$*$=\neg X^{\bullet}(t_1,...,t_n)$  ;  \\
- If $A=B \rightarrow C$, then $A$*$=B$*$ \rightarrow C$* ; \\
- If $A=\forall xB$ (resp. $A=\forall XB$), then $A$*=$\forall xB$* (resp. $A$*$=\forall XB$*) ;\\
- If $A=\forall X_C B$, then $A$*=$\forall X^{\bullet} B$*.\\
 
We have the following result (Nour 1997a). 

\begin{theorem}
 Let $A$ be a $\forall$-positive type of $AF2$ and $t$ a $\beta$-normal $\lambda 
C$-term. \\ If $\vdash_{M2} t:A$, then $t$ is a normal $\lambda$-term, and $\vdash_{AF2} t:A$. 
\end{theorem}

With each  predicate variable $X$ of $C2$, we associate a classical variable
$X_C$ having the same arity as $X$. For each formula $A$ of $C2$, we define the formula $A^C$ of
$M2$ in the following way : \\
- If $A=D(t_1,...,t_n)$ where $D$ is a constant symbol, then $A^C=A$ ;  \\
- If $A=X(t_1,...,t_n)$ where $X$ is a predicate symbol, then $A^C=X_C(t_1,...,t_n)$  ; \\  
- If $A=B \rightarrow C$, then $A^C=B^C \rightarrow C^C$ ; \\
- If $A=\forall xB$, then $A^C=\forall xB^C$ ;\\
- If $A=\forall XB$, then $A^C=\forall X_CB^C$.\\

As for relation betwen the systems $C2$ and $M2$, we have (Nour 1997a) :

\begin{theorem} Let $A$ be a type of $C2$, and $t$ a $\lambda C$-term. $\vdash_{C2} t:A$ iff
$\vdash_{M2} t:A^C$.  
\end{theorem}

\subsection{The integers in $M2$}

According to the results of the subsection 8.1, we obtain some results concerning 
integers in system $M2$ (Nour 1997a).

\begin{theorem}
Let $n\in {\bf N}$, if $\vdash_{M2} t :N[s^n (0)]$, then, $t \simeq_{\beta} \underline{n}$. 
\end{theorem}

Let $n\in {\bf N}$. By Theorem 8.1, a  classical integer of value $n$ is a
closed $\lambda C$-term $\theta_n$ such that $\vdash_{M2} \theta_n :N^C[s^n(0)]$. For the classical integers we
have only one operational characterization. In order to give this characterization, we shall need
some definitions.\\

Let $V$ be the set of variables of $\lambda C$-calculus. Let $P$ be an infinite
set of constants called stack constants \footnote{The notion of stack constants is taken from
a manuscript of Krivine.}. We define a set of $\lambda C$-terms $\Lambda CP$ by :\\
 - If $x \in V$, then $x \in \Lambda CP$ ; \\
 - If $t \in \Lambda CP$, and $x \in V$, then $\lambda  xt \in \Lambda CP$ ; \\
 - If $t \in \Lambda CP$, and $u \in \Lambda CP \bigcup P$, then $(t)u \in \Lambda CP$. \\
In other words, $t \in \Lambda CP$ iff the stack constants are in argument positions in
$t$. \\
We consider, on the set $\Lambda CP$, the following rules of reduction :\\
(1) $(\lambda  xu)tt_1...t_n \rightarrow (u[t/x])t_1...t_n$ for all $u,t \in \Lambda CP$ and
$t_1,...,t_n \in \lambda  CP \bigcup P$ ; \\
(2) $(C)tt_1...t_n \rightarrow (t)\lambda  x(x)t_1...t_n$ for all $t \in \Lambda  CP$ and $t_1,...,t_n \in
\lambda  CP \bigcup P$, and $x$  being $\lambda$-variable not appearing in $t_1,...,t_n$.\\ 
For any  $t,t' \in \Lambda CP$, we shall write $t \rhd_C t'$, if $t'$ is obtained from $t$ by
applying these rules finitely many times.\\

Let $\theta_1=\lambda x \lambda  f(C)\lambda  y(y)(f)(C)\lambda  z(y)(f)x$. We have $\vdash_{M2} \theta_1: N^C[s0]$ \\ 
$(\theta_1)xgp_0 \rhd_C (g)t_1 p_0$ ; $(t_1)p_1 \rhd_C (g)t_{2}p_0$ and $(t_2)p_2 \rhd_C (x)p_2$. In
general we have the following result (Nour 1997a). 

\begin{theorem}
Let $n \in {\bf N}$, $\theta_n$ a classical integer of value $n$, and
$x,g$ two distinct variables. \\
- If $n=0$, then, for every stack constant $p$, we have : $(\theta_n)xgp \rhd_C
(x)p$. \\
- If $n \not = 0$, then there is an $m \in {\bf N}$*, and a mapping $IÊ:Ê\{0,...,m\}\rightarrow {\bf N}$,
such that for all distinct stack constants  $p_0,p_1,...,p_m$, we have :\\
$(\theta_n)xgp_0 \rhd_C (g)t_1 p_{r_0}$ ;
$(t_i)p_i \rhd_C (g)t_{i+1}p_{r_i}$ $(1 \leq i \leq m-1)$ ; 
$(t_m)p_m \rhd_C (x)p_{r_m}$ 
where $I(0)=n$, $I(r_m)=0$, and $I(i+1)=I(r_i)-1$ $(0 \leq i \leq m-1)$.  
\end{theorem}

Theorem 8.4 allows to find the value of a classical integer. Let $\theta_n$ be a 
classical integer of value $n$. Let $p$ be a stack constant and $g,x$ two distinct variables. If
$(\theta_n)xgp \rhd_C (x)p$, then $n=0$. If not there is an $m \in {\bf N}$*, a sequence $(r_i)_{1\leq i
\leq m}$ where $(0 \leq r_i \leq m)$ and a mapping $JÊ:Ê\{0,...,m\}\rightarrow {\bf N}$ such that $J(0)=0$, and
$J(i+1)=J(r_i)+1$ $(0 \leq i \leq m-1)$. Therefore $J(r_m)=n$. 

\subsection{Storage operators for classical integers}

In system $M2$ we have a similar result to Theorem 5.2 (Nour 1997a).\\
 
Let $T$ be a closed $\lambda C$-term. We say that $T$ is a  storage operator for
classical integers iff for every $n \in {\bf N}$, there is a $\lambda C$-term $\tau_n
\simeq_{\beta} \underline{n}$, such that for every  classical integers $\theta_n$ of value $n$, there is a
substitution $\sigma$, such that $(T)\theta_n f \succ_C (f)\sigma(\tau_n)$ (where $f$ is a new variable).

\begin{theorem} 
If $\vdash_{M2} T: \forall x \{ N^C[x] \rightarrow \neg\neg N[x] \}$, then $T$ is a storage operator for classical
integers. 
\end{theorem}

Theorem 8.5 means that if $\vdash_{M2} T: \forall x \{ N^C[x] \rightarrow \neg\neg N[x] \}$, then $T$ 
takes a classical integer as an argument and return the Church integer corresponding to its value. It
is enough to do the proof of this Theorem in the propositionnal case.ÊThe type system $M$ is the
subsystem of $M2$ where we only have propositional variables and constants. We write $\Gamma\vdash_M
t:A$ if $t$ is typable in $M$ of type $A$ in the context $\Gamma$. Let $N = \forall X \{X,(X
\rightarrow X) \rightarrow X \}$. Theorem  8.5 is a consequence of the following Theorem (Nour 1997a).

\begin{theorem} 
If $\vdash_M T: N^C \rightarrow \neg\neg N$, then for every $n \in {\bf N}$, there is an $m \in {\bf N}$ and a $\lambda 
C$-term $\tau_m \simeq_{\beta} \underline{m}$, such that for every classical integer $\theta_n$ of value $n$,
there is a substitution $\sigma$, such that $(T)\theta_n f \succ_C (f)\sigma(\tau_m)$.  
\end{theorem}

Indeed, if $\vdash_{M2} T: \forall x \{ N^C[x] \rightarrow \neg\neg N[x] \}$, then $\vdash_M T:
N^C \rightarrow \neg\neg N$. Therefore for every $n \in {\bf N}$, there is an $m \in {\bf N}$ and $\tau_m \simeq_{\beta}
\underline{m}$, such that for every classical integer $\theta_n$ of value $n$, there is a substitution $\sigma$,
such that $(T)\theta_n f \succ_C (f)\sigma(\tau_m)$. We have $\vdash_{M2} \underline{n} : N^C[s^n(0)]$, then $ f:\neg
N[s^n(0)] \vdash_{M2} (T) \underline{n} f :\perp$, therefore $ f:\neg N[s^n(0)]\vdash_{M2} (f)\underline{m} :\perp$ and
$\vdash_{M2} \underline{m} : N[s^n(0)]$. Therefore $n=m$, and $T$ is a storage operator for classical
integers. \\

The proof of Theorem 8.6 uses two independent Theorems : the first one (Theorem 8.4)
expresses a property of classical integers and the second one (Theorem 8.7) expresses a property of
a $\lambda C$-terms of type $N^C \rightarrow \neg\neg N$.\\

Let $\nu$ and $f$ be two fixed variables. We denote by $x_{n,a,b,\overline{c}}$ (where $n$ is an integer,
$a,b$ two $\lambda$-terms, and $\overline{c}$ a finite sequence of $\lambda$-terms) a variable which does not
appear in $a,b,\overline{c}$. We have (Nour 1997a) :

\begin{theorem}
Let $\vdash_M T: N^C \rightarrow \neg\neg N$ and $n \in {\bf N}$. There is $m \in {\bf N}$ and a finite sequence of
head reductions $\{ U_i \succ_C V_i \}_{1\leq i\leq r}$ such that :\\  
1) $U_1 = (T)\nu f$ and $V_r = (f)\tau_m$ where $\tau_m \simeq_{\beta}\underline{m}$ ;\\
2) $V_i = (\nu) a b \overline{c}$ or $V_i = (x_{l,a,b,\overline{c}}) \overline{d}$  $(0 \leq l \leq
n-1)$ ;\\ 
3) If $V_i = (\nu)a b \overline{c}$, then $U_{i+1} = (a)\overline{c}$ if $n=0$ and
$U_{i+1} = ((b)x_{n-1,a,b,\overline{c}})\overline{c}$ if $n \neq 0$ \\ 
4) If $V_i = (x_{l,a,b,\overline{c}})\overline{d}$ $(0 \leq l \leq n-1)$, then $U_{i+1} = (a)\overline{d}$ if $l=0$
and $U_{i+1} = ((b)x_{l-1,a,b,\overline{d}})\overline{d}$ if $l \neq 0$.
\end{theorem}

Let $T$ be a closed $\lambda C$-term, and $D,E$ two closed types of $AF2$ type system.
We say that $T$ is a storage operator for the pair of types $(D,E)$ iff for every $\lambda$-term
$\vdash_{AF2} t:D$, there is $\lambda$-term $\tau'_t$ and $\lambda C$-term $\tau_t$, such that $\tau'_t
\simeq_{\beta} \tau_t$, $\vdash_{AF2} \tau'_t:E$, and for every $\vdash_{C2} \theta_t:D$, there is a substitution
$\sigma$, such that $(T)\theta_t f \succ_C (f)\sigma(\tau_t)$ (where $f$ is a new variable).\\

We can generalize Theorem 8.5 (Nour 1997a).

\begin{theorem}
Let $D,E$ two $\forall$-positive closed types of $AF2$ type system, such that $E$ does not contain
$\perp$. If $\vdash_{M2} T: D^C \rightarrow \neg\neg E$, then $T$ is a storage operator for the pair $(D,E)$.   
\end{theorem}

\section{The $\lambda \mu$-calculus}

\subsection{Pure and typed $\lambda \mu$-calculus}

$\lambda \mu$-calculus has two distinct alphabets of variables : the set of $\lambda$-variables
$x,y,z,...$, and the set of $\mu$-variables $\alpha,\beta,\gamma$,.... Terms (also called $\lambda
\mu$-terms) are defined by the following grammar : $t :=$ $x$ $\mid$ $\lambda xt$ $\mid$ $(t)t$
$\mid$  $\mu\alpha[\beta]t$.\\ 
The reduction relation of $\lambda \mu$-calculus is induced by fives
different notions of reduction :\\  
{\bf The computation rules \/}\\
($C_1$) $(\lambda xu)v \rightarrow u[v/x]$ \\
($C_2$) $(\mu \alpha u)v \rightarrow \mu \alpha u[v/$*$\alpha]$  where $u[\overline{v}/$*$\alpha]$ is obtained from $u$ by
replacing inductively each subterm of the form $[\alpha]w$ by $[\alpha](w)\overline{v}$.  \\
{\bf The simplification rules \/}\\
($S_1$)  $[\alpha]\mu \beta u \rightarrow  u[\alpha/\beta]$  \\
($S_2$)  $\mu \alpha [\alpha]u \rightarrow u$, if $\alpha$ has no free occurence in $u$  \\ 
($S_3$)  $\mu \alpha u \rightarrow \lambda x \mu \alpha u[x/$*$\alpha]$, if $u$ contains a subterm of the
form $[\alpha]\lambda yw$.  \\

Parigot (1992) has shown that :

\begin{theorem}
In $\lambda \mu$-calculus, reduction is confluent.
\end{theorem} 

The notation $u \succ_{\mu} v$ means that $v$ is obtained from $u$ by some head
reductions. The head equivalence relation is denoted by : $u \sim_{\mu} v$ iff there
is a $w$,  such that $u \succ_{\mu} w$ and $v \succ_{\mu} w$.\\

Proofs are written in a natural deduction system with several conclusions, presented
with sequents. One deals with sequents such that :\\  
- Formulas to the left of $\vdash$ are labelled with $\lambda$-variables ;\\  
- Formulas to the right of $\vdash$ are labelled with $\mu$-variables, except one formula
which is labelled with a $\lambda \mu$-term ;\\  
- Distinct formulas never have the same label.\\
  
Let $t$ be a $\lambda \mu$-term, $A$ a type, $\Gamma = x_1:A_1,...,x_n:A_n$,
and $\Delta = \alpha_1:B_1,...,\alpha_m:B_m$. We define by means of the following
rules the notion ``$t$ is of type $A$ in $\Gamma$ and $\Delta$''. This notion is denoted by
$\Gamma\vdash_{FD2} t:A,\Delta$. \\
The rules (1),...,(8) of $AF2$ type system and the following rule :
\begin{center}
$  ( 9 ) \quad \displaystyle\frac{ \Gamma\vdash_{FD2} t:A,\beta:B,\Delta} { \Gamma\vdash_{FD2}\mu\beta
[\alpha]t:B,\alpha:A,\Delta}$ \\ \end{center}
Weakenings are included in the rules (2) and (9).\\

As in typed $\lambda$-calculus on can define $\neg A$ as $A  \rightarrow \perp$ and use the previous rules with
the following special interpretation of naming for $\perp$ : for $\alpha$ a $\mu$-variable, $\alpha :
\perp$ is not mentioned. This typed $\lambda$-calculus system is called $FD2$. It has the following
properties (Parigot 1992).

\begin{theorem}
1) Type is preserved during reduction.\\
2) Typable $\lambda \mu$-terms are strongly normalizable. 
\end{theorem}

\subsection{Classical integers}

Let $n$ be an integer. A classical integer of value $n$ is a closed
$\lambda \mu$-term $\theta_n$ such that $\vdash_{FD2} \theta_n :N[s^n(0)]$.  \\
Let $x$ and $f$ fixed variables, and $N_{x,f}$ be the set of $\lambda \mu$-terms defined by the
following grammar : $u$ $:=$ $x$ $\mid$ $(f)u$ $\mid$ $\mu\alpha[\beta]x$ $\mid$
$\mu\alpha[\beta]u$.\\ 
We define, for each $u \in N\sb{x,f}$ the set $rep(u)$, which is intuitively
the set of integers potentially repesented by $u$ :\\
- $rep(x) = \{ 0 \}$ ;\\
- $rep((f)u) = \{ n+1$ if $n \in rep(u) \}$ ;\\
- $rep(\mu\alpha[\beta]u)= \bigcap rep(v)$ for each subterm $[\alpha]v$ of $[\beta]u$.\\

The following Theorem characterizes the classical integers (Parigot 1992).

\begin{theorem}
The normal classical integers of value $n$ are the $\lambda \mu$-terms of the form $\lambda x\lambda f
u$ with $u \in N_{x,f}$ without free $\mu$-variable and such that $rep(u)=\{n\}$.
\end{theorem}

Let $\theta = \lambda x\lambda f u$ where\\
$u=(f)\mu\alpha[\alpha](f)\mu\phi[\alpha](f)\mu\psi[\alpha](f)(f)
\mu\beta[\phi](f)\mu\delta[\beta](f)\mu\gamma[\alpha](f)\mu\rho[\beta](f)x$.\\ 
We can check that $rep(u)=\{4\}$. Then $\theta$ is a classical integer of value 4.\\

We will present now a simple method to find the value of a classical integer. We define, for each
$u\in N\sb{x,f}$ the set $val(u)$, which is intuitively the set of the possible values of $u$ : \\
- $val(x)=\{ 0 \}$ ; \\
- $val((f)u)=\{ n+1$ if $n\in val(u) \}$ ;\\
- $val(\mu\alpha[\beta]u)=\bigcup val(v)$ for each subterm $[\alpha]v$ of $[\beta]u$.\\\
Let $u\in N\sb{x,f}$ without free $\mu$-variable and $\alpha_1,...,\alpha_n$ the
$\mu$-variables of $u$ which satisfy :  $\alpha_1$ is the $\mu$-variable such that
$[\alpha_1](f)^{i_1 }x$ is a subterm of $u$, $\alpha_j$ $(2\leq j\leq n)$ is the $\mu$-variable such
that $[\alpha_j](f)^{i_j}\mu\alpha_{j-1}u_{j-1}$ is a subterm of $u$, and $u=(f)^{i_{n+1}}\mu\alpha_n
u_n$. Let $t_0 =x$ and $t_j=\mu\alpha_j u_j$ $(1\leq j \leq n)$. \\

We have (Nour 1997b).

\begin{lemma}
For every $(1\leq j \leq n+1)$ :\\  
1) $val(t_{j-1})=\displaystyle{\{\sum_{1 \leq k \leq j} \: i_k \}}$.\\ 
2) For each subterm $t$ of $u_j$, such that $t \neq (f) ^r t_k$  $(0\leq k\leq j-1)$, 
$val(t)=\emptyset$. \\
In particular $val(u)= \displaystyle{\{ \sum_{1 \leq k \leq n+1} \: i_k \}} $.
\end{lemma}

Using Lemma 9.1 and the fact that for each $u \in N_{x,f}$, $rep(u) \subseteq val(u)$ we deduce the
following result (Nour 1997b) :

\begin{theorem}
If $\theta$ is a normal classical integer of value $n$, then $\theta=\lambda x\lambda fu$ with $u \in N_{x,f}$
without free $\mu$-variable and such that $val(u)=\{ n \}$. 
\end{theorem}

Then to find the value of a normal classical integer $\theta = \lambda x\lambda f u$, we try the
$\mu$-variables $\alpha_j$ $(1\leq j\leq n+1)$ and the integers $i_j$ $(1\leq j\leq n+1)$ of the
$\lambda \mu$-term $u$. The value of $\theta$ is equal to $\displaystyle{ \sum_{1 \leq k \leq n+1} \:
i_k }$.

\subsection{Storage operators in $\lambda \mu$-calculus}

Let $T$ be a closed $\lambda$-term. We say that $T$ is a storage operator for
classical integers iff for every $(n \geq 0)$, there is $\lambda$-term $\tau_n
\simeq_{\beta} \underline{n}$, such that for every  classical integers $\theta_n$ of value $n$, there is a
substitution $\sigma$, such that $(T)\theta_n f  \sim_{\mu} \mu \alpha [\alpha] (f)\sigma(\tau_n)$
(where $f$ is a new variable).\\

Parigot (1993a) has shown that :

\begin{theorem}
If $\vdash_{AF2} T: \forall x \{ N^g[x] \rightarrow \neg\neg N[x] \}$, then $T$ is a storage operator for
classical integers. 
\end{theorem}

In order to define, in this framework, the equivalent of system $M2$, the
demonstration of $\neg \neg A \rightarrow A$ should not be allowed for all formulas $A$, and
thus we should prevent the occurrence of some formulas on the right. Thus the following
definition. \\ 

We add to the $FD2$ type system the new following rules :\\

\begin{minipage}[t]{170pt}
$  ( 6' ) \quad \displaystyle\frac{ \Gamma\vdash t:A, \Delta } { \Gamma\vdash t: \forall X_C A, \Delta }$  {\rm (*)}\\
\end{minipage} 
\begin{minipage}[t]{170pt}\sl 
$  ( 7') \quad  \displaystyle\frac{ \Gamma\vdash t: \forall X_C A, \Delta} { \Gamma\vdash t:A[G/ X_C], \Delta }$  {\rm
(**)}\\  \end{minipage}
With the following conditions : (*) $X_C$ has no free occurence in $\Gamma$ and (**) $G$ is a classical
type.\\
We call $M2$ the new type system, and we write $\Gamma \vdash_{M2} t:A , \Delta$ if $t$ is of type $A$
in the $\Gamma$ and $\Delta$.\\
Let $T$ be a closed $\lambda \mu$-term. We say that $T$ is a storage operator for classical
integers iff for every $(n \geq 0)$, there is $\lambda \mu$-term $\tau_n \simeq_{\beta}
\underline{n}$, such that for every  classical integers $\theta_n$ of value $n$, there is a substitution
$\sigma$, such that $(T)\theta_n f  \sim_{\mu} \mu \alpha [\alpha] (f)\sigma(\tau_n)$ (where $f$ is a new variable).\\

We have the following result :

\begin{theorem} 
If $\vdash_{M2} T: \forall x \{ N^C[x] \rightarrow \neg\neg N[x] \}$, then $T$ is a storage operator for classical
integers.  \end{theorem}

\end{document}